\documentclass[12pt]{amsart}
\usepackage{latexsym}
\usepackage{amsmath}
\usepackage{amsfonts}
\usepackage{amssymb}
\usepackage{amsthm}
\theoremstyle{plain}
\newtheorem{theorem}{Theorem}[section]
\newtheorem{lemma}[theorem]{Lemma}
\newtheorem{proposition}[theorem]{Proposition}

\newtheorem*{theorem*}{Theorem}

\theoremstyle{definition}

\newtheorem*{example}{Example}
\newtheorem*{acknowledgment}{Acknowledgment}

\theoremstyle{remark}

\newtheorem*{notation}{Notation}

\newcommand{\E}{\mathbb{E}}

\newcommand{\N}{\mathbb{N}}

\newcommand{\cC}{\mathcal{C}}

\newcommand{\cG}{\mathcal{G}}
\newcommand{\cI}{\mathcal{I}}

\newcommand{\cX}{\mathcal{X}}

\newcommand{\cY}{\mathcal{Y}}
\newcommand{\cZ}{\mathcal{Z}}

\newcommand{\nnorm}[1]{\lvert\!|\!| #1|\!|\!\rvert}
\newcommand{\type}[1]{^{[#1]}}
\newcommand{\bx}{\mathbf{x}}

\begin{document}
\title{Convergence of multiple
ergodic averages for some commuting transformations}

\author{Nikos Frantzikinakis and Bryna Kra}

\address{Department of Mathematics, McAllister Building,
The Pennsylvania State University,
University Park, PA 16802}
\email{nikos@math.psu.edu}
\email{kra@math.psu.edu}
\date{}

\begin{abstract}
We prove the $L^{2}$ convergence for the linear multiple ergodic
averages of commuting transformations $T_{1}, \ldots, T_{l}$,
assuming that each map $T_i$ and each pair $T_iT_j^{-1}$ is
ergodic for $i\neq j$. The limiting behavior of  such averages is
controlled by a particular factor, which is an inverse limit of
nilsystems. As a corollary we show that the limiting behavior of
linear multiple ergodic averages is the same for  commuting
transformations.
\end{abstract}

\maketitle

\section{Introduction}
We consider the multiple ergodic averages
\begin{equation}\label{eq:main}
\frac{1}{N}\sum_{n=0}^{N-1}f_1(T_1^nx)\cdot
f_2(T_2^nx)\cdot\ldots\cdot f_l(T_l^nx) \ ,
\end{equation}
where $T_1, T_2, \ldots, T_l$ are commuting measure preserving
transformations of a probability space $(X, \mathcal{X}, \mu)$.
Such averages, with $T_1 = T, T_2 = T^2, \ldots, T_l=T^l$, were
originally studied by Furstenberg~\cite{F2} in his proof of
Szemer\'edi's Theorem, and for general commuting transformations
by Furstenberg and Katznelson~\cite{FK} in their proof of the
multidimensional Szemer\'edi Theorem. The convergence in
$L^2(\mu)$ for the first case was proved in~\cite{HK}.

With certain hypotheses on the transformations, convergence for
three commuting transformations was proven by Zhang~\cite{zhang}.
Under the same hypotheses, we generalize this convergence
result for $l$ commuting transformations:

\begin{theorem}\label{T:main}
Let $l\geq 1$ be an integer. Assume that $T_1, T_2, \ldots, T_l$
are commuting invertible ergodic measure preserving
transformations of a measure space $(X, \cX, \mu)$ so that
$T_iT_j^{-1}$ is ergodic for all $i,j\in \{1, 2, \ldots, l\}$ with
$i\neq j$.  Then if  $f_1, f_2, \ldots, f_l\in L^\infty(\mu)$ the
averages
$$
 \frac{1}{N}\sum_{n=0}^{N-1}f_1(T_1^nx)\cdot
f_2(T_2^nx)\cdot\ldots \cdot f_l(T_l^nx)
$$
converge in $L^2(\mu)$ as $N\to +\infty$.
\end{theorem}

In order to prove convergence, we show that in the average
\eqref{eq:main} we can replace each function by its conditional
expectation on a certain factor and that this factor is an inverse
limit of translations on a nilmanifold.  However, this does not
hold for general commuting transformations. In
Section~\ref{S:cfactors} we show that certain ergodic assumptions
on the differences $T_iT_j^{-1}$, such as the ones we impose in
Theorem~\ref{T:main}, are necessary for the characteristic factor
to be an inverse limit of transformations on a nilmanifold. The
general case for commuting transformations remains open and only
is known for $l=2$~\cite{CL}. Convergence for certain distal
systems was obtained in~\cite{Les}.

Using the machinery developed for the proof of Theorem
\ref{T:main} and a result of Ziegler~\cite{ziegler}, in
Section~\ref{sec:commutingsystems} we prove an identity
illustrating that for ergodic commuting transformations the
limiting behavior of the corresponding linear multiple ergodic
averages is the same:
\begin{theorem}\label{T:2}
Suppose that   $T$ and $S$ are commuting invertible measure
preserving transformations of a probability space $(X, \cX, \mu)$
and that both $T$ and $S$ are ergodic. Then for all integers
$l\geq 1$ and any $f_1,\ldots,f_l\in L^\infty(\mu)$ we have:
\begin{equation}\label{E:same}
\lim_{N\to\infty}\frac{1}{N}\sum_{n=0}^{N-1}f_1(T^nx)\cdot \ldots
\cdot f_l(T^{ln}x)=
\lim_{N\to\infty}\frac{1}{N}\sum_{n=0}^{N-1}f_1(S^nx)\cdot \ldots
\cdot f_l(S^{ln}x)
\end{equation}
in $L^2(\mu)$.
\end{theorem}

We summarize the strategy of the proof of Theorem~\ref{T:main}.
The factors that control the limiting behavior of the
averages~\eqref{eq:main} are the same as those that arise in the
proof of the convergence of multiple ergodic averages along
arithmetic progressions in~\cite{HK}. These factors are defined by
certain seminorms and the starting point for the convergence of
the averages~\eqref{eq:main} is  that commuting ergodic systems
have the same associated seminorms (Proposition~\ref{L:same}). We
use this observation in Proposition \ref{P:bounds} to bound the
limsup of the $L^2$-norm of the averages by the seminorms  of the
individual functions. The structure theorem of \cite{HK} then
implies that it suffices to check convergence when every
transformation $T_i$ is isomorphic to an inverse limit of
translations on a nilmanifold. The main technical difficulty is to
prove that this isomorphism can be taken to be simultaneous and
that the factor sub-$\sigma$-algebras that appear in the inverse
limits can be chosen to be the same for all  transformations
$T_i$ (Theorems \ref{T:inverse} and \ref{T:nilrotation}).

\begin{notation}
For $l$ commuting maps $T_1, \ldots, T_l$, we use the shorthand
$(X,\mu, T_1, \ldots, T_l)$, omitting the $\sigma$-algebra, to denote the measure preserving
system obtained by the maps $T_1, \ldots, T_l$ acting on a fixed
measure space $(X, \cX, \mu)$. For simplicity of notation, we
assume that all functions in this article are real valued.  With
minor modifications, the same definitions and results hold for
complex valued functions.
\end{notation}

\section{Seminorms and factors $Z_k$}
Assume that $(X,\mu,T)$ is an ergodic system; we summarize a
construction and some results from \cite{HK}.

For an integer $k\geq 0$, we write $X\type k=X^{2^k}$ and $T\type k:
X\type k\to X\type k$ for the map $T\times T\times\ldots\times T$, taken
$2^k$ times.  We let $V_k = \{0,1\}^k$ and write elements of $V_k$
without
commas or parentheses.  Elements of $X^{[k]}$ are written $\bx =
(x_{\epsilon}:\epsilon\in V_k)$.

We define a probability measure $\mu\type k$ on $X\type k$, that
is invariant under $T\type k$, by induction. Set $\mu\type 0=\mu$.
For $k\geq 0$, let $\cI\type {k}$ be the $\sigma$-algebra of
$T\type {k}$-invariant subsets of $X\type {k}$. Then $\mu\type
{k+1}$ is the relatively independent square of $\mu\type{k}$ over
$\cI\type{k}$. This means that if $F',F''$ are bounded functions
on $X\type{k}$ then
\begin{equation}
\label{eq:def-muk} \int_{X\type
{k+1}}F'(\bx')F''(\bx'')\,d\mu\type {k+1}(\bx', \bx''):=
\int_{X\type{k}}\E(F'\mid\cI\type{k})\,\E(F''\mid\cI\type{k})
\,d\mu\type{k}\ .
\end{equation}

For a bounded function $f$ on $X$ and integer $k\geq 1$ we define
\begin{equation}
\label{eq:def-seminorm} \nnorm f_{k}^{2^{k}}= \int_{X\type
{k}}\prod_{j=0}^{2^{k}-1} f(x_j)\, d\mu\type {k}(\bx) \ ,
\end{equation}
and we note that the integral on the right hand side is
nonnegative. We immediately have that $\nnorm f_1=|\int f \
d\mu|$.  It is shown in~\cite{HK} that for every integer $k\geq
1$, $\nnorm\cdot_k$ is a seminorm on $L^\infty(\mu)$ and
using the Ergodic Theorem it is easy to check that
\begin{equation}
\label{eq:recur} \nnorm f_{k+1}^{2^{k+1}} =\lim_{N\to+\infty}\frac
1N\sum_{n=0}^{N-1} \nnorm{f\cdot T^nf}_{k}^{2^{k}}\ .
\end{equation}
Furthermore, it is shown in \cite{HK} that for every integer
$k\geq 1$ the seminorms define factors $Z_{k-1}(X)$  in the
following manner: the $T$-invariant sub-$\sigma$-algebra
$\cZ_{k-1}(X)$ is characterized by
\begin{equation}\label{eq:def-Z_k}
\text{ for } f\in L^\infty(\mu),\  \E(f|\cZ_{k-1})=0\text{ if and
only if } \nnorm f_{k} = 0 \ ,
\end{equation}
then  $Z_{k-1}(X)$ is defined to be the factor of $X$ associated
to the sub-$\sigma$-algebra $\cZ_{k-1}(X)$. Thus defined, $Z_0(X)$
is the trivial factor, $Z_1(X)$ is the Kronecker factor and more
generally, $Z_k(X)$ is a compact abelian group extension of
$Z_{k-1}(X)$. We denote the restriction of $\mu$ to $\cZ_{k}(X)$
by $\mu_k$.

Define $\cG = \cG(X,T)$ to be the group of measure preserving
transformations $S$ of $X$ which satisfy for every integer $k\geq
1 $: the transformation $S^{[k]}$ of $X^{[k]}$ leaves the measure
$\mu^{[k]}$ invariant and acts trivially on the invariant
$\sigma$-algebra $\cI^{[k]}$.

In our context, there is more than one transformation acting on
the space and so we need  some notation to identify the
transformation with respect to which the seminorm is defined. We
write $\nnorm \cdot_{k,T}$ for the \emph{$k$-th seminorm} with
respect to the map $T$.  Similarly we write $Z_{k}(X,T)$ for the
\emph{$k$-th factor} associated to $T$, $\mu^{[k]}_T$ for the
measure defined by  \eqref{eq:def-muk}, and $\cI^{[k]}(T)$ for the
$T^{[k]}$-invariant subsets of $X^{[k]}$.

\section{Characteristic factors for commuting
transformations}\label{S:cfactors} Let  $T_1, T_2, \ldots, T_l$ be
commuting invertible measure preserving transformations of a
probability space $(X, \cX, \mu)$. We say that a
sub-$\sigma$-algebra $\cC$ of $\cX$ is a \emph{characteristic
factor for $L^2(\mu)$-convergence} of the averages \eqref{eq:main}
if $\cC$ is $T_i$ invariant for $i=1,\ldots,l$, and the averages
\eqref{eq:main} converge in $L^2(\mu)$ to zero whenever
$\mathbb{E}(f_i|\cC)=0$ for some $i\in\{1,\ldots,l\}$. In this
section we show that under the assumptions of
Theorem~\ref{T:main}, the sub-$\sigma$-algebra $\cZ_{l-1}$ is
characteristic for $L^2(\mu)$-convergence of the above averages.
Before proving this
 we give an example showing that this  does not hold for general
commuting transformations.

\begin{example} Consider the measure preserving transformations
$T_1=R\times S_1$, $T_2=R\times S_2$, acting on the measure space
$(X=Y\times Z,\mathcal{Y}\times \mathcal{Z},\mu\times\nu)$ and
suppose
 that $T_1$ and $ T_2$   commute. Taking $l=2$,
$f_1(y,z)=f(y)$, and $f_2(y,z)=\overline{f}(y)$  we see that the
limit of the averages in \eqref{eq:main} is zero if and only if
$f$ is zero almost everywhere. It follows that any
sub-$\sigma$-algebra that is characteristic  for
$L^2(\mu)$-convergence of these averages contains $\cY$. So if
$\cY\not\subset \cZ_1(Y,R)$, meaning that if $(Y,\mu,R)$ is not a
Kronecker system for ergodic $R$, then  $\cZ_1(X,T_1)$ is not a
characteristic factor. More generally, if there exists a function
that is $T_1T_2^{-1}$ invariant but not
$\mathcal{Z}_1(X,T_i)$-measurable, for $i=1,2$, the same argument
gives that $\cZ_1(X,T_i)$ is not a characteristic factor  for
$i=1,2$.
\end{example}

We turn now to the proof of the main result of the section. The
most important ingredient of the proof is the next observation.
This  result, as well as  Proposition~\ref{P:bounds}, were also
obtained independently by I.~Assani (personal communication).

\begin{proposition}\label{L:same}
 Assume that $T$ and $S$ are commuting  measure
preserving transformations of a  probability space $(X, \cX, \mu)$
and that both $T$ and $S$ are ergodic. Then $S\in \cG(X,T)$,
$\cG(X,T)=\cG(X,S)$, and for all integers $k\geq 1$ and all $f\in
L^{\infty}(\mu)$,  $\mu^{[k]}_T=\mu^{[k]}_S$,
$\cI^{[k]}(T)=\cI^{[k]}(S)$, $\nnorm f_{k,T} = \nnorm f_{k, S}$,
and $ Z_k(X,T)=Z_k(X,S)$.
\end{proposition}
\begin{proof}
By Lemma $5.5$ in~\cite{HK}, $S\in \cG(X,T)$.   We use induction
on $k$ to show that $\mu^{[k]}_T=\mu^{[k]}_S$ and
$\cI^{[k]}(T)=\cI^{[k]}(S)$. The statement is obvious for $k=0$.
Suppose that it holds for some  integer $k\geq 1$. By
Equation~\eqref{eq:def-muk}, we have that
$\mu^{[k+1]}_T=\mu^{[k+1]}_S$.
 Since $S\in \cG(X,T)$, we have that
 $S^{[k+1]}$ leaves the
measure $\mu^{[k+1]}_T=\mu^{[k+1]}_S$ invariant and acts
trivially on $\mathcal{I}_T^{[k+1]}$. Hence, $\cI^{[k+1]}(T)
\subset \cI^{[k+1]}(S)$. Reversing the roles of $T$ and $S$, we have
that $\cI^{[k+1]}(T)=\cI^{[k+1]}(S)$. This completes the
induction.

By the definition of $\cG(X,T)$ and Equations~\eqref{eq:def-seminorm}
and~\eqref{eq:def-Z_k}, the equalities
 $\cG(X,T)=\cG(X,S)$, $\nnorm f_{k,T} =
\nnorm f_{k, S}$  and  $Z_k(X,T)=Z_k(X,S)$
follow.
\end{proof}
\begin{proposition}\label{P:bounds}
Let $l\geq 1$ be an integer. Assume that $T_1, T_2, \ldots, T_l$
are commuting invertible ergodic measure preserving
transformations of a probability space $(X, \cX, \mu)$ such that
$T_iT_j^{-1}$ is ergodic for all $i,j\in \{1, 2, \ldots, l\}$ with
$i\neq j$.

(i) If $f_1, f_2, \ldots, f_l\in L^\infty(\mu)$, then
$$
\limsup_{N\to +\infty}
\Big\|\frac{1}{N}\sum_{n=0}^{N-1}f_1(T_1^nx)\cdot
f_2(T_2^nx)\cdot\ldots\cdot f_l(T_l^nx)\Big\|_{L^2(\mu)} \leq
\min_{i=1, 2, \ldots, l}\nnorm {f_i}_{l} \ ,
$$
where the seminorm is taken with respect to any $T_i$, $i\in\{1,
2, \ldots, l\}$.

(ii) The sub-$\sigma$-algebra $\cZ_{l-1}$ is a characteristic
factor for convergence of the averages \eqref{eq:main}.
\end{proposition}

\begin{proof}
Note that by Proposition~\ref{L:same}, we can take these seminorms
and the sub-$\sigma$-algebras $\cZ_l$ with respect to any of the
maps $T_i$ with $i\in\{1, 2, \ldots, l\}$, since they are all the
same.

Part $(ii)$ follows immediately from part $(i)$ and the definition
of the sub-$\sigma$-algebra $\cZ_{l-1}$. We prove the inequality
of part $(i)$ by induction. For $l=1$, this follows from the
ergodic theorem.

Assume the statement holds for $l-1$ functions and that $f_1,f_2,
\ldots, f_l\in L^\infty(\mu)$ with $\|f_j\|_\infty \leq 1$ for $j
= 1, 2, \ldots, l$. First assume that $i\in\{2, 3, \ldots, l\}$.
(The case $i=1$ is similar, with the roles of $T_1$ and $T_2$
reversed.) Let
$$
u_n = f_1(T_1^nx)\cdot f_2(T_2^nx)\cdot\ldots\cdot f_l(T_l^nx) \ .
$$
For convenience we use the notation $T^nf(x)=f(T^nx)$. By the Van
der Corput Lemma (see Bergelson~\cite{bergelson}),
\begin{equation}\label{E:u_n}
 \limsup_{N\to\infty}\Bigl\lVert \frac{1}{N}\sum_{n=0}^{N-1}
u_{n}\Bigr\lVert_{L^2(\mu)}^2\leq \limsup_{M\to
\infty}\frac{1}{M}\sum_{m=0}^{M-1} \limsup_{N\to\infty}\Big|
\frac{1}{N}\sum_{n=0}^{N-1} \langle{ u_{n+m}, u_n \rangle  }\Big|\
.
\end{equation}
A simple computation gives that
$$
\frac{1}{N}\sum_{n=0}^{N-1} \langle{ u_{n+m}, u_n \rangle }
\leq  \Bigl\lVert\frac{1}{N}\sum_{n=0}^{N-1}(T_2T_1^{-1})^n (
T_2^mf_2 f_2)\cdot \ldots\cdot (T_lT_1^{-1})^n( T_l^mf_l f_l
)\Bigr\lVert_{L^2(\mu)}.
$$
For $j=2, 3, \ldots, l$, define $S_j = T_jT_1^{-1}$.  The maps
$S_2, S_3, \ldots, S_l$ commute, since the maps $T_1, T_2, \ldots,
T_l$ do, and they also commute with $T_1, \ldots, T_l$.
Furthermore for $i\neq j$, the transformation $S_iS_j^{-1} =
T_iT_j^{-1}$ is ergodic by assumption.  Hence, we can use the
inductive assumption and \eqref{eq:recur} to bound the right hand
side in \eqref{E:u_n} by
$$
\limsup_{M\to \infty} \frac{1}{M}\sum_{m=0}^{M-1}\nnorm
{f_iT_i^mf_i}_{l-1, S_i} \leq \limsup_{M\to \infty}
\Big(\frac{1}{M} \sum_{m=0}^{M-1}\nnorm {f_iT_i^mf_i}_{l-1,
S_i}^{2^{l-1}}\Big)^{1/2^{l-1}} =\nnorm {f_i}_{l, T_i}^2\ ,
$$
where the last equality follows from~\eqref{eq:recur} and
Proposition~\ref{L:same}.
\end{proof}


\section{Structure of the characteristic factor and the proof of convergence}
\label{sec:commutingsystems} We have shown that under the
assumptions of Theorem~\ref{T:main}, a characteristic factor for
$L^2(\mu)$-convergence of the averages \eqref{eq:main} is
$Z_{l-1}$. The structure theorem of \cite{HK} gives  that for
$i=1,\ldots,l$ the system $(Z_{l-1}, \mu_{l-1}, T_i)$ is
isomorphic to an inverse limit of ($l-1$)-step nilsystems. But
this isomorphism apriori depends on the transformation $T_i$, and
the sub-$\sigma$-algebras that appear in the inverse limits depend
on the transformations. In this and the following section we
extend several results from \cite{HK} and use them to deal with
these technical difficulties.

Throughout this section, $T_1,\ldots,T_l$ are commuting, ergodic,
measure preserving transformations of a probability space
$(X,\mathcal{X},\mu)$. If $S_1,\ldots,S_l$ are measure preserving
transformations of a probability space $(Y,\mathcal{Y},\nu)$, we
say that the system $(X,\mu, T_1, \ldots, T_l)$ is
\emph{isomorphic} to $(Y,\nu, S_1, \ldots, S_l)$ if there exist
sets $X'\subset X$, $Y'\subset Y$ of full measure that are
invariant under all the transformations on their respective
spaces, and  a measurable bijection $\phi\colon X' \to Y'$
carrying $\mu$ to $\nu$ that satisfies $\phi(T_i(x))=S_i(\phi(x))$
for  every $x\in X'$ and all $i=1,\ldots,l$.  When $\phi$ is not
assumed to be injective, then we say that $(Y, \nu, S_1, \ldots,
S_l)$ is a \emph{factor} of $(X,\mu, T_1, \ldots, T_l)$. For
simplicity of notation we assume that $X=X'$ and $Y=Y'$.

We say that the system $(X,\mu,T_1, \ldots, T_l)$ is an
\emph{extension} of its factor $(Y,\nu,S_1, \ldots, S_l)$ (let
$\pi\colon X\to Y$ denote the factor map) \emph{by a compact
abelian group} $(V,+)$  if there exist measurable cocycles
$\rho_1,\ldots,\rho_l\colon Y \to V$ and a measure preserving
bijection $\phi\colon X \to Y \times V$ (we let $m_V$ denote the Haar measure
on $V$), satisfying:

\noindent  $(i)$ $\phi$ preserves $Y$, meaning that $\phi^{-1}(\cY
\times V)=\pi^{-1}(\cY)$ up to sets of measure zero, where
$\cY\times V =\{A\times V\colon A\in \cY\}$, and

\noindent $(ii)$ $\phi( T_i(x))= T_i'(\phi(x))$ for all $x\in X$,
$y\in Y$, $v\in V$ and $i\in\{1,\ldots,l\}$, where
$$T_i'(y,v)=(S_i(y),v+\rho_i(y))\ .$$
Setting  $\tilde{\rho}=(\rho_1,\ldots,\rho_l)$, we let
$Y\times_{\tilde{\rho}} V$ denote the system $(Y\times V,\nu\times
m_V,T_1',\ldots,T_l')$.


We say that the system $(X,\mu, T_1, \ldots, T_l)$ has
\emph{order} $k$ if $X = Z_k(X)$. (Note that $Z_k(X) = Z_k(X,T_i)$
does not depend on $i$ by Proposition~\ref{L:same}.) It has
\emph{toral Kronecker factor} if $(Z_1,\mu_1,T_1,\ldots, T_l)$ is
isomorphic to a system $(V,m_V,R_1, \ldots, R_l)$, where $V$ is a
compact abelian Lie group, $m_V$ is the Haar measure, and $R_i$,
for $i=1, \ldots, l$, is a rotation on $V$. Finally, it is
\emph{toral} if  it is of order $k$ for some integer $k\geq 1$, it
has toral Kronecker factor, and for $ j=1,\ldots, k-1$, the system
$(Z_{j+1},\mu_{j+1},T_1,\ldots, T_l)$ is an extension of
$(Z_j,\mu_j,T_1, \ldots, T_l)$ by a finite dimensional torus.

We say that the system $(X,\mu,T_1,\ldots, T_l)$ is an
\emph{inverse limit} of a sequence of factors
$\{(X_j,\mu_j,T_1,\ldots,T_l)\}_{j\in \N}$, if
$\{\mathcal{X}_j\}_{j\in\N}$ is an increasing sequence of
sub-$\sigma$-algebras  invariant under the transformations
$T_1,\ldots, T_l$, and  such that
$\bigvee_{j\in\N}\mathcal{X}_j=\mathcal{X}$ up to sets of measure
zero. If in addition for every $j\in\N$ the factor system
$(X_j,\mu_j,T_1,\ldots, T_l)$ is isomorphic to a toral system of
order $k$, we say that $(X,\mu,T_1,\ldots, T_l)$ is an
\emph{inverse limit of a sequence of toral systems of order} $k$.

The next result extends Theorem 10.3  from \cite{HK} to the case
of several commuting transformations and we postpone the proof
until Section~\ref{S:inverse}.
\begin{theorem}\label{T:inverse}
Any system $(X,\mu,T_1,\ldots, T_l)$ of order $k$ is an inverse
limit of a sequence
$\{(X_i,\mu_i,T_1,\ldots,T_l)\}_{i\in\mathbb{N}}$ of toral systems
of order $k$.
\end{theorem}
To describe the characteristic factors, we briefly review the
definition of a nilsystem. Given a group $G$, we denote the
commutator of $g,h\in G$ by $[g,h]=g^{-1}h^{-1}gh$. If $A,
B\subset G$, then $[A,B]$ is defined to be the subgroup generated
by the set of commutators  $\{[a,b]:a\in A, b\in B\}$.  Set
$G^{(1)} = G$ and for integers $k\geq 1$, we inductively define
$G^{(k+1)}= [G, G^{(k)}]$.
 A group $G$
is said to be {\em $k$-step nilpotent} if its $(k+1)$ commutator
$[G, G^{(k)}]$ is trivial.  If $G$ is a $k$-step nilpotent Lie
group and $\Gamma$ is a discrete cocompact subgroup, then the
compact space $X = G/\Gamma$ is said to be a {\em $k$-step
nilmanifold}. The group $G$ acts on $G/\Gamma$ by left translation
and the translation by a fixed element $a\in G$ is given by
$T_{a}(g\Gamma) = (ag) \Gamma$.  There exists a unique probability
measure $m_{G/\Gamma}$ on $X$ that is invariant under the action
of $G$ by left translations (called the {\em Haar measure}).
Fixing elements $a_1,\ldots, a_l \in G$, we call the system
$(G/\Gamma,m_{G/\Gamma}, T_{a_1},\ldots, T_{a_l})$ a {\em $k$-step
nilsystem} and call each map $T_a$ a {\em nilrotation}.

The next result extends the structure theorem of \cite{HK}
(Theorem 10.5) to the case of several ergodic commuting
transformations:
\begin{theorem}\label{T:nilrotation}
 Let $l\geq 1$ be an integer and let $(X,\mu, T_1,\dots,T_l)$ be a toral
 system of order $k$. Then there exist a $k$-step nilpotent Lie group
 $G$, a discrete and cocompact subgroup $\Gamma$ of $G$, and commuting
 elements $a_1,\dots,a_l$ of $G$, such that the system $(X,\mu,
 T_1,\dots,T_l)$ is isomorphic to the $k$-step nilsystem
 $(G/\Gamma, m_{G/\Gamma},T_{a_1},\dots,T_{a_l})$.
\end{theorem}

\begin{proof}
Let $\mathcal{G}=\mathcal{G}(X,T_1)$. By the structure theorem
in~\cite{HK} the group $\mathcal{G}$ is $k$-step nilpotent, and by
Proposition \ref{L:same} we have that $T_1,\ldots,T_l
\in\mathcal{G}$. Let $G$ be the subgroup of $\mathcal{G}$ that is
spanned   by the connected component of the identity and the
transformations $T_1,\ldots,T_l$. At the end of the proof of
Theorem~10.5 in~\cite{HK}, it is shown that there exist a discrete
cocompact subgroup $\Gamma$ of $G$ and a measurable bijection
$\phi\colon G/\Gamma\to X$ that carries the Haar measure
$m_{G/\Gamma}$ to $\mu$, such that for all $S\in G$ the
transformation $\phi^{-1} S \phi$ of $G/\Gamma$ is the left
translation on $G/\Gamma$ by $S$. Since $\phi$ does not depend on
$S\in G$, and $T_1,\ldots,T_l \in G$, the proof is complete.

\end{proof}

We now combine the previous results to prove Theorem~\ref{T:main}:

\begin{proof}[Proof of Theorem~\ref{T:main}]
If $\mathbb{E}(f_i|\cZ_{l-1})=0$  for some $i\in\{1,\ldots, l\}$,
then $\nnorm{f_i}_l = 0$ and  by Proposition~\ref{P:bounds} the
limit of the averages \eqref{eq:main} is zero. Hence,  we can
assume that $X=Z_{l-1}$. By Theorem~\ref{T:inverse}, the system
$(X,\mu,T_1,\ldots,T_l)$ is an inverse limit of toral systems
$(X_i,\mu_i,T_1,\ldots,T_l)$ of order $l-1$. It follows from
Theorem \ref{T:nilrotation} that $(X_i,\mu_i,T_1,\ldots,T_l)$ is
isomorphic to an $(l-1)$-step nilsystem. Using an approximation
argument, we can assume that $X=G/\Gamma$, $\mu=m_{G/\Gamma}$, and
the transformations $T_i$ are given by nilrotations $T_{a_i}$ on
$G/\Gamma$ where $a_i$ are commuting elements of $G$ for
$i=1,\ldots, l$. If $a=(a_1,\ldots,a_l)\in (G/\Gamma)^l$, then by
\cite{Lei} the average
$$
\lim_{N\to \infty} \frac{1}{N}\sum_{n=0}^{N-1} T_a^n f
$$
converges everywhere in $(G/\Gamma)^l$ as $N\to\infty$ for $f$
continuous. Using the convergence on the diagonal and a standard
approximation argument the result follows.
\end{proof}

We also now have the tools to prove Theorem \ref{T:2}:

\begin{proof}[Proof of Theorem~\ref{T:2}]
The $L^2(\mu)$ convergence for these averages was proved in
\cite{HK}. From Proposition~\ref{L:same} we have that
$Z_{l-1}(X,T)=Z_{l-1}(X,S)$,
 and  by Proposition~\ref{P:bounds}  this common factor  is
 characteristic for $L^2$-convergence of the averages in
\eqref{E:same}. Hence, we can assume that $X=Z_{l-1}$. Then by
Theorem \ref{T:inverse}, the system $(X,\mu,T,S)$ is an inverse
limit of toral systems $(X_i,\mu_i,T,S)$ of order $l-1$. It
follows by Theorem \ref{T:nilrotation} that $(X_i,\mu_i,T,S)$ is
isomorphic to an $(l-1)$-step nilsystem. Using an approximation
argument, we can assume $X=G/\Gamma$, $\mu=m_{G/\Gamma}$, and the
transformations $T$ and $S$ are given by ergodic nilrotations
$T_a,T_b$ on $G/\Gamma$.

In~\cite{ziegler}, Ziegler gives a formula for the limit of the
averages in~\eqref{E:same} when the transformations are
nilrotations. Using this identity it is clear that the limit is
the same for all ergodic nilrotations on $G/\Gamma$. Since both
$T_a$ and $T_b$ are ergodic, the result follows.
\end{proof}

\section{Proof of Theorem~\ref{T:inverse}}\label{S:inverse}
The proof of Theorem~\ref{T:inverse}   is carried out in three
steps, as in~\cite{HK} for a single transformation.  We give all
the statements of the needed modifications, but only include the
proof when it is not a simple rephrasing of the corresponding
proof in ~\cite{HK}. We summarize the argument. Suppose that
the system $(X,\mu,T_1,\ldots,T_l)$ is of order $k+1$,  meaning
that $X=Z_{k+1}$. First we prove that the system
$(Z_{k+1},\mu_{k+1},T_1,\ldots,T_l)$ is a connected compact
abelian group extension of $(Z_{k},\mu_k,T_1,\ldots,T_l)$. Next we
show that if $(Z_{k},\mu_k,T_1,\ldots,T_l)$ is an inverse limit of
systems $\{(Z_{k,i},\mu_{k,i},T_1,\ldots,T_l)\}_{i\in
\mathbb{N}}$, then the cocycle that defines the extension  is
measurable with respect to some $\cZ_{k,i}$. Finally, we combine
the first two steps to complete the proof by induction. To prove
Lemmas~\ref{L:connected} and~\ref{L:torus}, we make use of the
analogous results in \cite{HK} for a single transformation,
extending them to several transformations. The argument given in
Step 3 is similar to the one used in \cite{HK} and we include it
as it ties together the previous lemmas.

{\bf Step 1.} We extend Theorem $9.5$ from \cite{HK}.
\begin{lemma} \label{L:connected}
Let $k,l\geq 1$ be integers and let $T_1,\ldots,T_l$ be commuting
ergodic measure preserving transformations  of a probability space
$(X,\mathcal{X},\mu)$. Then the system
$(Z_{k+1},\mu_{k+1},$ $T_1,\ldots,T_l)$ is an extension of the system
$(Z_{k},\mu_k,T_1,\ldots,T_l)$ by a connected compact abelian
group.
\end{lemma}
\begin{proof}
 Using the analogous result for a single transformation
 in~\cite{HK} (part $(ii)$ of Theorem $9.5$),  we get that
  there exist a
connected compact abelian group $V$, a cocycle $\rho_1\colon
Z_k\to V$, a measure preserving bijection $\phi\colon Z_{k+1}\to
Z_k\times V$ that preserves $Z_k$ such that $T_1(\phi(x))=
\phi(T_1'(y,v))$ for $y\in Z_k$, $v\in V$, and a measure
preserving transformation $S_1:V\to V$ such that
$$T_1'(y,v)=(S_1(y),v+\rho_1(y))\ .$$
For $i=2,\ldots,l$, define $T_i'=\phi^{-1} T_i \phi$. Since both
$\phi$ and $T_i$ preserve $Z_k$, we have that $T_i'$  has the form
$$
T_i'(y,v)=(S_i(y),Q_i(y,v))\ ,
$$
for some measurable transformation $Q_i:Z_k\times V\to  V$ for
$i=2,\ldots, l$. By Corollary $5.10$ in~\cite{HK} all maps
$R_u\colon Z_k\times V\to Z_k\to V$ defined by $R_u(y,v)=(y,v+u)$
 belong to the center of $\mathcal{G}(T_1')$. By
Proposition~\ref{L:same}, $T_i'\in \mathcal{G}(T_1')$ and so
$T_i'$ commutes with $R_u$ for all $u\in V$. This can only happen
if $Q_i$ has the form $Q_i(y,v)=v+\rho_i(y)$ for some measurable
$\rho_i\colon Y\to V$. This completes the proof.
\end{proof}

{\bf Step 2.} For $i=1,\ldots,l$ define the coboundary operator
$\partial_i$ by $\partial_i(f)=f\circ T_i-f$. If $U$ is a compact
abelian group and $\rho_i\colon X\to U$ are cocycles, then the
cocycle $\tilde{\rho}=(\rho_1,\ldots,\rho_l)$ is called an
   \emph{l-cocycle} if
   $$
   \partial_i\rho_j=\partial_j \rho_i , \text{ for all } \
   i,j\in \{1,\ldots,l\} \ .
   $$
    This is equivalent to
  saying that the maps $S_1,\ldots,S_l$  defined on $X\times U$
  by  $S_i(x,u)=(T_ix,u+\rho_i(x))$ commute.

If $U$ is a compact abelian group and $\rho_1, \ldots, \rho_l:X\to
U$ are cocycles, then the cocycle $\tilde{\rho} = (\rho_1, \ldots,
\rho_l)$ is an $\emph{l-coboundary}$ for the system
$(X,\mu,T_1,\ldots,T_l)$  if there exists a measurable function
$f\colon X\to U$ such that $\rho_i=\partial_i f$ for
$i=1,\ldots,l$. Furthermore, $\tilde{\rho}$ is said to be an
$\emph{l-quasi coboundary}$ for $(X,\mu,T_1,\ldots,T_l)$ if there
exists a measurable function $f\colon X\to U$ and $c_i\in U$ such
that $\rho_i=c_i+\partial_i f$ for all $i=1,\ldots,l$.  Two
$l$-cocycles are \emph{cohomologous} if their difference is an
$l$-coboundary.

For clarity of exposition, we include a few more definitions
from~\cite{HK}. Let $(X,\mu,T)$ be a measure preserving system and
let $U$ be a compact abelian group. We say that the cocycle
$\rho\colon X\to U$ is  \emph{ergodic} with respect to the system
$(X,\mu,T)$ if the extension $(X\times U, \mu\times m_U,
T_{\rho})$, where $T_{\rho}:X\times U\to X\times U$ is given by
$T_{\rho}(x,u)= (Tx, u+\rho(x))$, and $m_U$ is the Haar measure on
$U$, is ergodic.
 An $l$-cocycle
$\tilde{\rho}=(\rho_1,\ldots,\rho_l)$ is ergodic
 with respect to the system $(X,\mu,T_1,\ldots,T_l)$ if $\rho_i$ is
 ergodic with respect to $(X,\mu,T_i)$ for $i=1,\ldots,l$.

For an integer $k\geq 1$ and $\epsilon\in V_k$, write $|\epsilon|
= \epsilon_1+\ldots+\epsilon_k$ and $s(\epsilon)
=(-1)^{|\epsilon|}$. For each $k\geq 1$, we define the map
$\Delta^k\rho:X^{[k]}\to U$ by
$$\Delta^k\rho(\bx) = \sum_{\epsilon\in
V_k}s(\epsilon)\rho(x_{\epsilon}) \ .
$$
We say that the cocycle $\rho\colon X\to U$ is \emph{of type $k$}
with respect  to the system $(X,\mu,T)$ if the cocycle $\Delta
^k\rho:X^{[k]}\to U$ is a coboundary of the system $(X^{[k]},
\mu^{[k]}, T^{[k]})$. An $l$-cocycle $\tilde{\rho}\colon X\to U^l$
is said to be of \emph{type $k$} with respect to the system
$(X,\mu,T_1,\ldots,T_l)$ if each coordinate cocycle is of type
$k$.

Let $(Y,\mathcal{Y},\nu)$ be a probability space, let $V$ be a
compact abelian group with Haar measure $m_V$ and let $X=Y\times
V$. The action $\{ R_v: v\in V\}$  of measure preserving
transformations  $R_v\colon X\to X$ defined by $R_v(y,u)=(y,u+v)$
is called an action on $X$ by \emph{vertical rotations} over $Y$.

Let $(X,\mathcal{X},\mu)$ be a probability space, let $V$ be a
connected compact abelian group and let $\{S_v:v\in V\}$ be an
action of $V$ on $X$ by measure preserving transformations
$S_v\colon X \to X$. The action $S_v$ is said to be \emph{free} if
there exists a probability space
 $(Y,\mathcal{Y},\nu)$ and an action $R_v\colon Y\times V\to Y\times V$  by
 vertical rotations over $Y$
 such that the actions $\{S_v:v\in V\}$ and $\{R_v:v\in V\}$ are isomorphic. This means that there
  exists a measurable
 bijection $\phi:Y\times V \to X$, mapping $\nu\times m_V$ to
 $\mu$ and satisfying  $\phi( R_v(y,u))=S_v(\phi(y,u))$ for all $y\in Y$ and
all $u,v\in V$.

\begin{lemma}\label{L:torus}
Let $l\geq 1$ be an integer and let $T_1,\ldots,T_l$ be commuting
ergodic measure preserving transformations  of a  probability
space $(X,\mathcal{X},\mu)$.  Let $\{S_v: v\in V\}$ be a free
action of a compact abelian group $V$ on $X$ that commutes with
$T_i$ for $i=1, \ldots, l$. Let $U$ be a finite dimensional torus
and let $\tilde{\rho}=(\rho_1,\ldots,\rho_l)\colon X\to U^l$ be an
$l$-cocycle of type $k$ for some integer $k\geq 2$. Then there
exists a closed subgroup $V'$ of $V$ such that $V/V'$ is a compact
abelian Lie group, and there exists an $l$-cocycle
$\rho'=(\rho_1',\ldots,\rho_l')$, cohomologous to $\rho$, such
that $\rho_i'\circ S_v=\rho_i'$ for every $v\in V'$.
\end{lemma}
\begin{proof}  As before, we define the
operators $\partial_i$ by $\partial_i(f)=f\circ T_i-f$ for
$i=1,\ldots,l$ and define $\partial_v(f)=f\circ S_v-f$ for $v\in
V$.  Using the analogous result for a single transformation in
\cite{HK} (Corollary 9.7), we get that for $i=1, \ldots, l$ there
exist closed subgroups $V_i$ of $V$ such that $V/V_i$ is a compact
abelian Lie group, and measurable $f_i\colon X\to U$ such that
$\bar{\rho}_i=\rho_i+\partial_if_i$ satisfies $\partial_v
\bar{\rho_i}=0$ for all $v\in V_i$. If $W=\bigcap_{i=1}^k V_i$,
then $V/W$ is a compact abelian Lie group and  all the previous
relations hold for $v\in W$. We take $V'$ to be   the
connected component of the identity element in $W$. For
$i=1,\ldots,l$ let $\rho_i'=\rho_i+\partial_if_1$. Since $V/W$ is
a compact abelian Lie group and $W/V'$ is finite, we have that
$V/V'$ is also a compact abelian Lie group. It suffices to show
that  for $v\in V'$ we have $\partial_v\rho_i'=0$ for all
$i=1,\ldots,l$.

Since $\partial_1\rho_i=\partial_i\rho_1$ and the operators
$\partial_i$ commute, it follows that
$\partial_1\rho_i'=\partial_i\rho_1'$. Moreover, since $S_v$
commutes with all the $T_i$, the operators $\partial_v$ and
$\partial_i$ commute. It follows that
$$
\partial_1\partial_v\rho_i'=\partial_v\partial_1\rho_i'=
\partial_v\partial_i\rho_1'= \partial_i\partial_v\rho_1'=0
$$
 for $v\in V'$, where
the last equality holds since by assumption, $\partial_v\rho_1'=0$
for $v\in V'$. Since $\partial_1\partial_v\rho_i'=0$ and $T_1$ is
ergodic, we have that for $i=1,\ldots,l$ and $v\in V'$ there
exists a constant $c_{v,i}\in U$ such that
$\partial_v\rho_i'=c_{v,i}$.  It follows that for fixed $i$ the
map $c_{v,i}\colon V'\to U$ is a measurable homomorphism.
Moreover, for $i=1,\ldots,l$ we have
$\rho_i'=\bar{\rho_i}+\partial_i(f_1-f_i)$, and $\partial_v
\bar{\rho_i}=0$ for all $v\in V'$.  Hence,
$c_{v,i}=\partial_v\rho_i'=\partial_ig_i$ for all $v\in V'$, where
$g_i=\partial_v(f_1-f_i)$. The spectrum of $T_i$ is countable and
so the last relation implies that $c_{v,i}$ can only take on
countably many values for $v\in V'$.  So  $c_{v,i}\colon V'\to U$
is a measurable homomorphism and takes on countably many values.
Since $V'$ is connected, it follows that $c_{v,i}=0$ for
$i=1,\ldots,l$ and $v\in V'$.
\end{proof}

Using the previous Lemma,
 the proof of the next result is identical to that of Lemma $10.4$
in~\cite{HK} and so  we  do not reproduce it.

\begin{lemma}\label{L:inverse}
Let $l\geq 1$ be an integer, let $T_1,\ldots,T_l$ be commuting
ergodic measure preserving transformations  of a  probability
space $(X,\mathcal{X},\mu)$, let $U$ be a finite dimensional
torus, and let $\tilde{\rho}\colon X\to U^l$ be an ergodic
$l$-cocycle of type $k$ and measurable with respect to $\cZ_k$ for
some integer $k\geq 1$. Assume that the system
$(X,\mu,T_1,\ldots,T_l)$ is an inverse limit of the sequence of
systems $\{(X_i,\mu_i,T_1,\ldots,T_l)\}_{i\in\N}$. Then
$\tilde{\rho}$ is cohomologous to a cocycle $\tilde{\rho}'\colon
X\to U^l$, which is measurable with respect to $\mathcal{X}_i$ for
some $i$.
\end{lemma}

{\bf Step 3.} We complete the proof of Theorem~\ref{T:inverse} by
using induction on $k$. If $k=1$ we can assume that $(X,\mu,T_1)$
is an ergodic rotation on a compact abelian group $V$. Since $T_i$
commutes with $T_1$ for $i=1,\ldots l$, it follows that each $T_i$
is also a rotation on $V$. A compact abelian group is a Lie group
if and only if its dual is finitely generated. Hence, every
compact abelian group is an inverse limit of compact abelian Lie
groups and the result follows.

 Suppose that the result holds for some integer $k\geq 1$. Assume that
 $(X,\mu,T_1,\ldots,T_l)$
 is a system of order $k+1$. By
Lemma~\ref{L:connected}, we get that the system
$(X,\mu,T_1,\ldots,T_l)$ is an extension of
$(Z_k,\mu_k,T_1,\ldots,T_l)$ by a connected compact abelian group
$V$. Thus we can assume that $X=Z_k\times_{\tilde{\rho}} V$ where
$\tilde{\rho}=(\rho_1,\ldots,\rho_l)\colon Z_k\to V^l$ is the
$l$-cocycle defining the extension. The cocycle $\tilde{\rho}$  is
an ergodic $l$-cocycle since every system $(X,\mu,T_i)$ is.
Moreover, since $(X,\mu,T_i)$ and $(Z_k,\mu_k,T_i)$ are systems of
order $k+1$,  by Corollary 7.7 in \cite{HK} the cocycle $\rho_i$
is  of type $k+1$ for $i=1,\ldots l$. Hence, $\tilde{\rho}$ is an
$l$-cocycle of type $k+1$.

Since $V$ is a connected compact abelian group it can be written
as an inverse limit of finite dimensional tori $V_j$. Let
$\tilde{\rho_j}\colon Z_k\to V_j^l$ be the projection of
$\tilde{\rho}$ on the quotient $V_j$ of $V$. By the inductive
hypothesis, $(Z_k,\mu_k,T_1,\ldots,T_l)$ can be written as an
inverse limit of toral systems
$(Z_{k,i},\mu_{k,i},T_1,\ldots,T_l)$. By Lemma \ref{L:inverse},
for every integer $j$ there exist an integer $i_j$, and an
$l$-cocycle $\tilde{\rho}_i'\colon Y\to V_j^l$ cohomologous to
$\tilde{\rho}_i$ and measurable with respect to $\cZ_{k,i_j}$.
Without loss of generality, we can assume that the sequence $i_j$
is increasing. Let $X_j=Z_k\times V_j$. Then the system
$(X_j,\mu_j,T_1,\ldots,T_l)$ is isomorphic to the toral system
$Z_{k,i_j} \times_{\tilde{\rho}_j'}V_j$. Since
$(X,\mu,T_1,\ldots,T_l)$ is an inverse limit of the sequence
$\{(X_j,\mu_j,T_1,\ldots,T_l)\}_{j\in\N}$, the proof is complete.

\begin{acknowledgment} We thank B.~Host for helpful discussions
during the preparation of this article.
\end{acknowledgment}

\end{document}